# PISOT NUMBERS AND PRIMES

Andrei Vieru

**Abstract**

We define and study a transform whose iterates bring to the fore interesting relations between Pisot numbers and primes. Although the relations we describe are general, they take a particular form in the Pisot limit points.

We give three elegant formulae, which permit to locate on the whole semi-line all limit points that are not integer powers of other Pisot numbers.

## 1. INTRODUCTION AND NOTATIONS

### 1.1. PRELIMINARY DEFINITIONS

#### 1.1.1. PISOT NUMBERS

A Pisot–Vijayaraghavan number, also called simply a Pisot number, is a real algebraic integer greater than 1 such that all its Galois conjugates are less than 1 in absolute value.

#### 1.1.2. PALINDROMIC, ANTI-PALINDROMIC AND SEMI-PALINDROMIC POLYNOMIALS

A polynomial with real coefficients $P(x) = a_n x^n + a_{n-1} x^{n-1} + \ldots + a_1 x + a_0$ is called palindromic if, for any $i$, $a_i = a_{n-i}$

A polynomial with real coefficients $P(x) = a_n x^n + a_{n-1} x^{n-1} + \ldots + a_1 x + a_0$ is called anti-palindromic if, for any $i$, $a_i = -a_{n-i}$

It is called semi-palindromic if, for any even $i < n$, $a_i = -a_{n-i}$ while, for $i=n$ and for any odd $i$, $a_i = a_{n-i}$

### 1.1.3. RECIPROCAL, ANTI-RECIPROCAL AND SEMI-RECIPROCAL POLYNOMIALS

A polynomial with real coefficients $P(x) = a_n x^n + a_{n-1} x^{n-1} + \ldots + a_1 x + a_0$ is the reciprocal polynomial of $P^*(x) = b_n x^n + b_{n-1} x^{n-1} + \ldots + b_1 x + b_0$ if, for any $i$, $a_i = b_{n-i}$. The polynomials $P(x)$ and $P^*(x)$ are called each other's anti-reciprocals if, for any $i$, $a_i = -b_{n-i}$. The polynomials $P(x)$ and $P^*(x)$ are called each other's semi-reciprocals if, for any even $i < n$, $a_i = -b_{n-i}$, while for any odd $i$ $a_i = b_{n-i}$ and if $a_n = b_0$ and if $a_0 = b_n$.

### 1.2. NOTATIONS

Let $[x]$ designate the *nearest integer function* and let $\theta$ be any *Pisot number* which is not a limit point[1] of the closed set of Pisot numbers. To avoid confusion, the brackets [ ] will be used only on this purpose. We shall use brackets { } to designate only sequences, for example $\{[\theta^n]\}_{n \in N}$

We'll use the parenthesis signs ( ) only to indicate the order of elementary arithmetic operations and to fill in the argument of a function or a transform.

Despite some ambiguities, we'll designate by $\theta_k$ any Pisot number of degree $k$ (i.e. with minimal polynomial of degree $k$) that is not a limit point.

We'll drop the index $k$ whenever there is no need of keeping it.

In this article we'll study the iterates of a transform of the sequences $\{\theta^n\}_{n \in N}$, namely $\{\theta^n\}_{n \in N} \to \{\theta^n(\theta^n - [\theta^n])\}_{n \in N}$.

The sequence $\{\theta^n(\theta^n - [\theta^n])\}_{n \in N}$ will be considered as the first iterate of the transform, and we shall designate it by $\{I(\theta^n)\}_{n \in N}$. Following these notations, $\{[I(\theta^n)]\}_{n \in N}$ is an **integer sequence,** while $\{I(\theta^n)\}_{n \in N}$ is, at best, an **almost integer sequence**.)

---

[1] unless otherwise indicated

$\{I^2(\theta^n)\}_{n \in N} = \{\theta^n(\theta^n(\theta^n - [\theta^n]) - [\theta^n(\theta^n - [\theta^n])])\}_{n \in N}$ will be regarded as the second iterate of the same transform.

Generally, if $\{I^j(\theta^n)\}_{n \in N}$ is the *j*-th iterate of this transform, we'll define the *(j+1)-th* iterate of this transform by the formula:

$$\{I^{j+1}(\theta^n)\}_{n \in N} = \{\theta^n(I^j(\theta^n) - [I^j(\theta^n)])\}_{n \in N}$$

Of course, $\{\theta^n\}_{n \in N}$ will be considered as the order 0 iterate of the above-mentioned transform.

For a given $\theta$, we shall sometimes use $u_m^k$ to designate the integer $[I^k(\theta^m)]$.

## 2. THE FIVE MAIN CONJECTURES

### 2.1. *THE PRIME DIVISIBILITY CONJECTURE*

For any Pisot number $\theta_n$ (which is not a limit point), there are at least two non negative integers $k$ ($k \leq n - 2$) and $m(\theta_n)$ such as

for any prime $p \geq m(\theta_n)$, $[I^k(\theta_n^p)] \equiv 0 \pmod{p}$

### 2.2. *THE CONGRUENCE MODULO PRIME CONJECTURE*

For any Pisot number $\theta_n$ (which is not a limit point) and for any given non negative integer $k \leq n-1$, there is an integer $m(k, \theta_n)$ such as one the four following statements holds:

a) for any prime $p \geq m(k, \theta_n)$ $[I^k(\theta_n^p)] \equiv 0 \pmod{p}$

b) for any prime $p \geq m(k, \theta_n)$ $[I^k(\theta_n^p)] \equiv 1 \pmod{p}$

c) for any prime $p \geq m(k, \theta_n)$ $[I^k(\theta_n^p)] \equiv -1 \pmod{p}$

## 2.3. THE LINEAR RECURRENCE CONJECTURE

For any Pisot number $\theta_n$ (be it a limit point or not) and for any non negative integer $k$ ($k < n-1$), there is an integer $m(k,\theta_n)$ and $j$ signed integers $b_1, b_2,\ldots, b_j$ such as for any $l \geq m(k,\theta_n)$

$$\left[I^k(\theta_n^l)\right] = b_1\left[I^k(\theta_n^{l-1})\right] + b_2\left[I^k(\theta_n^{l-2})\right] + \ldots + b_j\left[I^k(\theta_n^{l-j})\right] = \sum_{i=1}^{j} b_i\left[I^k(\theta_n^{l-i})\right]$$

## 2.4. THE POLYNOMIAL COEFFICIENTS CONJECTURE

For any Pisot number $\theta$ of degree $n \geq 3$ (here $\theta$ is supposed to be neither a limit point nor the plastic constant), root of the polynomial $P(x) = a_n x^n + a_{n-1} x^{n-1} + \ldots + a_1 x + a_0$, there are two integers $j(P)$ and $k(P)$ such as

a)  for any $m \geq k(P)$,

$$\left[I^{n-2}(\theta^m)\right] = (-1)^{i(1)} a_1 \left[I^{n-2}(\theta^{m-1})\right] + (-1)^{i(2)} a_2 \left[I^{n-2}(\theta^{m-2})\right] + \ldots + (-1)^{i(n)} a_n \left[I^{n-2}(\theta^{m-n})\right] =$$

$$\sum_{j=1}^{n} (-1)^{i(j)} a_j \left[I^{n-2}(\theta^{m-j})\right] \qquad (1°)$$

where $i$ is a function $i:\{a_1, a_2,\ldots, a_n\} \to \{1, 0\}$

The function $i$ is defined as follows:

If $n$ is odd, then for any $k$ ($1 \leq k \leq n$)  $i(k) = 0$

If $n$ is even, then for any even $k$ ($1 \leq k \leq n$)  $i(k) = 0$

and for any odd $k$ ($1 \leq k \leq n$) $i(k) = 1$

For any given $\theta$ different from the plastic number **(1°)** may be also written as:

$$u_m^{n-2} = (-1)^{i(1)} a_1 u_{m-1}^{n-2} + (-1)^{i(2)} a_2 u_{m-2}^{n-2} + \ldots + (-1)^{i(n)} a_n u_{m-n}^{n-2} = \sum_{j=1}^{n} (-1)^{i(j)} a_j u_{m-j}^{n-2}$$

b) Let again $\theta$ be any Pisot number of degree $n \geq 2$ (supposed not to be a limit point), root of the polynomial $P(x) = a_n x^n + a_{n-1} x^{n-1} + \ldots + a_1 x + a_0$, there is an integer $k_0(P)$ such as for any $m \geq k_0(P)$ the following formula holds:

$$u_m^0 = -a_0 u_{m-n}^0 - a_1 u_{m-(n-1)}^0 - a_2 u_{m-(n-2)}^0 - \ldots - a_{n-1} u_{m-1}^0 = -\sum_{j=0}^{n-1} a_j u_{m-n+j}^0$$

## 2.5. CHARACTERISTIC POLYNOMIALS

More generally, it seams that, for even $k$, the characteristic polynomial of the linear recurrence sequence $\left\{\left[I^{\frac{k}{2}+1}\left(\theta_k^n\right)\right]\right\}_{n \geq m(k,\,\theta_k)}$ is anti-palindromic if $\theta_k$ is either not a limit point or a limit point of type $\alpha$ (see below, **3.1.,** *a*)). It is semi-palindromic if $\theta_k$ is a limit point of type $\beta$ (see below, **3.1.,** *b*)).

For even or odd $k$, the characteristic polynomials for the sequences

$$\left\{\left[I^m\left(\theta_k^n\right)\right]\right\}_{n \geq m(k,\,\theta_k)} \text{ and } \left\{\left[I^{k-m+2}\left(\theta_k^n\right)\right]\right\}_{n \geq m(k,\,\theta_k)}$$

are either reciprocal or anti-reciprocal or semi-reciprocal (depending on evenness or oddness of $k$ and $m$ (and of the type of the Pisot number: limit point or not, type of limit point, etc.).

## 2.6. *THE CONSTANT* 1 *CONJECTURE*

For any Pisot number $\theta$ of degree $n \geq 3$ (which is not a limit point), root of the polynomial $P(x) = a_n x^n + a_{n-1} x^{n-1} + \ldots + a_1 x + a_0$, there is an integer $j(P)$ such as for any $m \geq j(P)$ $\left|I^{n-1}\left(\theta^m\right)\right| = 1$

## 3. ON THE PISOT LIMIT POINTS SMALLER THAN 2

## 3.1. LOGARITHMIC EQUATIONS FOR PISOT LIMIT POINTS

It is well-known that the set of Pisot numbers between 1 and 2 have infinitely many limit points (as well as on the rest of the semi-line). It is also well-known that the sequence of these limit points is made of[2]:

*a*) the roots $\alpha_n$ of the polynomials $1-x+x^n(2-x)$ $(n \geq 1)$

---

[2] See also Mohamed Amara, *Ensembles fermés de nombres algébriques,* 1966, *Annales Scientifiques de l'École Normale Supérieure, 3ᵉ série, t3, pp. 230-260.* This article has been kindly sent me by Philippe Gille and Jean-Paul Allouche.

b) the roots $\beta_n$ of the polynomials $\dfrac{1 - x^{n+1}(2 - x)}{1 - x}$ $(n \geq 1)$

c) the number $\hat{o}'_2 = 1.90516616775...$, root of the polynomial $1 - 2x^2 - x^3 + x^4$ (and, as matter of fact, the square power of the second smallest Pisot number $1.380277...$, root of the polynomial $x^4 - x^3 - 1$; it is well-known that if $\theta$ is a Pisot number, then for any $m \geq 2$  $\theta^m$ is a limit point).

This sequence of limit points verifies the inequalities:

$$\alpha_1 = \beta_1 < \alpha_2 < \beta_2 < \alpha_3 < \hat{o}'_2 < \beta_3 < ... < \alpha_n < \beta_n < \alpha_{n+1} < ... < 2 \quad (\blacklozenge)$$

The limit points satisfy the following nice identities:

$$\forall n, \quad \dfrac{\text{Log}(2 - \beta_n)^{-1}}{\text{Log}(\beta_n)} = n + 1 \tag{I}$$

$$\forall n, \quad \dfrac{\text{Log}(2 - \alpha_n)^{-1}}{\text{Log}(\alpha_n)} - \dfrac{\text{Log}(\alpha_n - 1)^{-1}}{\text{Log}(\alpha_n)} = n \tag{II}$$

Besides, $\dfrac{\text{Log}(2 - \alpha_3)^{-1}}{\text{Log}(\alpha_3)} - \dfrac{\text{Log}(\alpha_3 - \alpha_1)^{-1}}{\text{Log}(\alpha_3)} = 1$

$$\dfrac{\text{Log}(2 - \alpha_2)^{-1}}{\text{Log}(\alpha_2)} = \dfrac{5}{2} \qquad \dfrac{\text{Log}(\alpha_2 - 1)^{-1}}{\text{Log}(\alpha_2)} = \dfrac{1}{2}$$

and

$$\dfrac{\text{Log}(2 - \hat{o}'_2)^{-1}}{\text{Log}(\hat{o}'_2)} - \dfrac{\text{Log}(\hat{o}'_2 - 1)^{-1}}{\text{Log}(\hat{o}'_2)} = \dfrac{7}{2}$$

### 3.2. CONGRUENCE MODULO PRIME CONJECTURES FOR PISOT LIMIT POINTS SMALLER THAN 2

While reading the statements of this chapter, one has to keep in mind that the Pisot limit points $\alpha_n$ and $\beta_n$ are Pisot numbers of degree $n+1$.

**1°a)** For any positive integer $n \geq 2$, there is an integer $m(n)$ such as

for any prime $p \geq m(n)$

$$\left[I^{0}\left(\alpha_{n}^{p}\right)\right] \equiv 2 \pmod{p}$$

**1°b)**

for any prime $p \geq m(n)$  $\left[I^{n-1}(\alpha_{n}^{p})\right] \equiv -1 \pmod{p}$

**1°c)** For any prime $p \geq m(n)$, for any $n$ and

for any strictly positive $k < n-1$, $\left[I^{k}(\alpha_{n}^{p})\right] \equiv 0 \pmod{p}$

**1°d)** For any $n > 0$ and for any integer $l \geq m(n)$

$$\left[I^{2n}\left(\alpha_{2n}^{l}\right)\right] = 1$$

**1°e)** For any $n > 0$ and for any integer $l \geq m(n)$

$$\left[I^{2n-1}\left(\alpha_{2n-1}^{2l-1}\right)\right] = 1, \text{ while } \left[I^{2n-1}\left(\alpha_{2n-1}^{2l}\right)\right] = -1$$

(note that the left part of the last three equalities might be written without the brackets [ ] of the *nearest integer function*.)

**2°a)** For any positive integer $n > 1$ there is an integer $m(n)$ such as for any prime $p \geq m(n)$ and for any integer $k < n$

$$\left[I^{k}(\beta_{n}^{p})\right] \equiv 1 \pmod{p}$$

**2°b)** For any integer $n > 0$ and for any integer $l \geq m(n)$

$$\left[I^{2n}(\beta_{2n}^{l})\right] = 1$$

**2°c)** For any integer $n > 0$ and for any integer $l \geq m(n)$

$$\left[I^{2n-1}(\beta_{2n-1}^{2l-1})\right] = 1 \text{ while } \left[I^{2n-1}(\beta_{2n-1}^{2l})\right] = -1$$

(note that the left part of this equality might be written without the brackets [ ] of the *nearest integer function*.)

**3°)** Remarkably, although it is a limit point (of degree 4), $\hat{o}_{2}{'}$ behaves in some sense like an ordinary Pisot number:

For any prime $p > 11$, $\left[I^{2}(\hat{O}_{2}^{'p})\right] \equiv 0 \pmod{p}$,

but, on the other hand, for any integer $m \geq 11$

$$\left[I^{3}(\hat{O}_{2}^{'m})\right] = -1$$

For sufficiently big $l$, $[I^{n-1}(\theta'^{l})] = -1$ (where $n$ is the degree of $\theta'$) is a typical equality for limit points $\theta'$ of the form $\theta^m$. (Note again that the left part of this equality might be written without the brackets [ ] of the *nearest integer function*.)

**4°)** It is very likely that the **Linear Recurrence Conjecture** holds as well for all Pisot numbers which happen to be limit points. In other words, for any $k$ smaller than $n-2$, all sequences $\{[I^k(\theta_n^l)]\}_{l \in N}$ are *additive* beginning with some term of rank $m(k, \theta_n)$.

**5°)** The **Polynomial Coefficients Conjecture** holds as well for all Pisot numbers that are also limit points.

**5α°)** For the numbers $\alpha_n$ the **Polynomial Coefficients Conjecture** may be formulated as follows:

For any $n \geq 2$, there is an integer $j(n)$ such as for any $m \geq j(n)$,

$$[I^{n-1}(\alpha_n^m)] = [I^{n-1}(\alpha_n^{m-(n+1)})] + (-2)^{n+1}[I^{n-1}(\alpha_n^{m-n})] + (-1)^n[I^{n-1}(\alpha_n^{m-1})] \qquad (°2)$$

**5β°)** For the numbers $\beta_n$ the **Polynomial Coefficients Conjecture** should be formulated as follows:

For any odd $n \geq 2$, there is an integer $j(n)$ such as for any $m \geq j(n)$, we have

$$[I^{n-1}(\beta_n^m)] = [I^{n-1}(\beta_n^{m-(n+1)})] + \sum_{i=1}^{(n+1)/2}[I^{n-1}(\beta_n^{m-(2i-1)})] - \sum_{i=1}^{(n-1)/2}[I^{n-1}(\beta_n^{m-2i})]$$

$$= [I^{n-1}(\beta_n^{m-(n+1)})] + (-1)^{i+1}\sum_{i=1}^{n}[I^{n-1}(\beta_n^{m-i})] \qquad (°3)$$

while, for any even $n \geq 2$, there is an integer $j(n)$ such as for any $m \geq j(n)$, we have

$$[I^{n-1}(\beta_n^m)] = [I^{n-1}(\beta_n^{m-(n+1)})] - [I^{n-1}(\beta_n^{m-n})] - [I^{n-1}(\beta_n^{m-(n-1)})] -$$

$$\ldots - \left[I^{n-1}(\beta_n^{m-2})\right] - \left[I^{n-1}(\beta_n^m)\right] = \left[I^{n-1}(\beta_n^{m-(n+1)})\right] - \sum_{j=1}^{n}\left[I^{n-1}(\beta_n^{m-j})\right] \qquad (4°)$$

## 4. QUADRATIC PISOT NUMBERS

Let $\theta$ be a quadratic Pisot number, root of the polynomial $x^2 + a_1 x + a_0$.

The sequence $\left\{\left[I(\theta^n)\right]\right\}_{n \in N}$ is always $\left\{-a_0^n\right\}_{n \in N}$

By the way, 1.5617520677202972947… (which is a root of the 'atypical' polynomial $x^6 - 2x^5 + x^4 - x^2 + x - 1$ and which is not a limit point) is the only known – at least to me – Pisot non quadratic number, such as for any even $l \geq 44$

$$I^{n-1}(\theta^l) = -1 \qquad (*)$$

while for any odd $l \geq 43$

$$I^{n-1}(\theta^l) = +1 \qquad (**)$$

## 5. THE "CONVERGENCE" SPEED CONJECTURE

**5.1.** How fast sequences generated by the iterates of the *I*-transform approach integers?

Let $\theta_n$ be any Pisot number of degree $n$ (limit point or not). Let $k_1(P)$ be the smallest of the numbers $k(P)$ – respectively $j_1(n)$ the smallest of the numbers $j(n)$ – for which the formulas **(1°)**, **(2°)**, **(3°)** or **(4°)** in the ***Polynomial Coefficients Conjectures*** (see pages 4 and 8) hold true.

For any integer $\lambda \leq n - 1$, and for any integers $\sigma$ and $\nu$

$$k_1(P) \leq \sigma < \nu < \infty \Rightarrow \left|I^\lambda(\theta^\nu) - \left[I^\lambda(\theta^\nu)\right]\right| < \left|I^\lambda(\theta^\sigma) - \left[I^\lambda(\theta^\sigma)\right]\right|$$

Following our notations, if $\theta$ is a limit point, we'll write:

For any integer $\lambda \leq n - 1$, and for any integers $\sigma$ and $\nu$

$$j_1(n) \leq \sigma < \nu < \infty \Rightarrow \left|I^\lambda(\theta^\nu) - \left[I^\lambda(\theta^\nu)\right]\right| < \left|I^\lambda(\theta^\sigma) - \left[I^\lambda(\theta^\sigma)\right]\right|$$

## 5.2. THE ERRATIC BEHAVIOR CONJECTURE

$k_1(P)$ and $j_1(n)$ may be arbitrarily big. In other words, the sequences of Pisot numbers and the iterates of the transform $I$ give rise to an arbitrarily long initial subsequence with 'erratic' behavior, before they become truly additive sequences.

## 6. NON INTEGER LIMIT POINTS GREATER THAN 2

The logarithmic equations (I) and (II) can be generalized in the following way (considering integers $m \geq 2$, $n \geq 1$ and $l < m$):

$$\frac{\text{Log}(m-x)^{-1}}{\text{Log}(x)} = n \qquad (\clubsuit)$$

$$\frac{\text{Log}(m-x)^{-1}}{\text{Log}(x)} - \frac{\text{Log}(x-m+l)^{-1}}{\text{Log}(x)} = n \qquad (\heartsuit)$$

Besides, there is an another logarithmic equation, namely

$$\frac{\text{Log}(x-m)^{-1}}{\text{Log}(x)} = n \qquad (\spadesuit)$$

whose real roots – along with the roots of the equations of types ($\clubsuit$) and ($\heartsuit$) – are always Pisot numbers.

Moreover, they always are Pisot *limit points*.

We conjecture that the set of solutions of the three families of equations ($\clubsuit$), ($\heartsuit$) and ($\spadesuit$) contain all Pisot *limit points* that are not integer powers of Pisot numbers.

It is not difficult to establish ($\blacklozenge$)-type inequalities (see page 8) for limit points on the semi-line that verify the equations ($\clubsuit$), ($\heartsuit$) and ($\spadesuit$). We shall not do it here, the reader can easily write them by himself. Let us only say that, setting $m = m_0$, all solutions fall in the open[3] interval $]m_0-1, m_0[$.

Let us also point out that, for $m = m_0$, $l$ has $m_0-1$ possible values. It means that in the interval $]m_0–1, m_0[$ we will find $m_0–1$ infinite sequences of Pisot limit points of the $\alpha$-type (see page 6, *a*) and (II)). Although the set of Pisot numbers is countable and nowhere dense, one can say, that the set of limit points in the interval $]m_0–1, m_0[$ is 'less dense' than in the interval $]m_0, m_0+1[$ (in the sense that the average interval between two consecutive Pisot limit points that belong to some family of solutions of the families of equations ♣, ♥ or ♠ (indexed on *l*, *m* and *n*) is smaller in $]m_0–1, m_0[$ is 'less dense' than in the interval $]m_0, m_0+1[$ .

---

[3] We do not consider here integer limit points : they are trivial.

For the Pisot limit point $\theta$ that satisfies the equation

$$\frac{\text{Log}(m_0 - x)^{-1}}{\text{Log}(x)} - \frac{\text{Log}(x - m_0 + l_0)^{-1}}{\text{Log}(x)} = n_0 \qquad (\heartsuit_0)$$

we'll state:

**1°°)** For any positive integer $n \geq 2$, there is an integer $m_0(n)$ such as

for any prime $p \geq m_0(n)$ $\left[I^{n-2}(\theta^p)\right] \equiv -1 \pmod{p}$

while,

for any $n \geq 2$ and for any integer $l \geq m_0(n)$

$$I^{n-1}(\theta^l) = 1$$

$$\left[I^0(\theta^p)\right] \equiv m_0 \pmod{p}$$

and, for any prime $p \geq m_0(n)$ and for any strictly positive $k$ smaller than $n - 2$

$$\left[I^k(\theta^p)\right] \equiv 0 \pmod{p}$$

Obviously, **1°°)** happens to be a generalization of **1°*a*)** (see end of page 6).

**Conclusion**

The methods described in this paper provide an algorithm able, for each given *n*, to generate integer sequences with the following property: For any prime *p*, $u_p \equiv n \pmod{p}$

This article describes how the iterates of the ***I***-transform generate integer sequences with remarkable properties. Almost nothing is said about the way the iterates of the ***I***-transform generate new recurrence formulae – those of the

(conjecturally) linear recurrence sequences that appear at each step of the iterating process.

Examples and results of computations may soon be found in the appendix of the extended version of this paper at http://www.andreivieru.com/page.php?id=4
Interesting examples of fixed points in the fractional part of numbers obtained through the iteration of the *I* transform arise in some computations. We believe infinitely many Pisot numbers yield such fixed points (for finitely many integer powers of them). Integer sequences that count subsets (with specific properties) of cyclic graphs appear in some computations.

Andrei Vieru
e-mail address: andreivieru007@gmail.com
www.andreivieru.com

**Acknowledgment**
We wish to express our gratitude to Vlad Vieru for computer programming.